\documentclass[11pt]{amsart}
\usepackage{graphicx}
\newtheorem{thm}{Theorem}[section]
\newtheorem{cor}[thm]{Corollary}
\newtheorem{lem}[thm]{Lemma}
\newtheorem{prop}[thm]{Proposition}
\theoremstyle{definition}

\theoremstyle{remark}

\numberwithin{equation}{section}

\newcommand{\shat}{{$(\hat{s})$}}
\newcommand{\stilde}{{$(\tilde{s})$}}
\newcommand{\s}{$(s)$}
\newcommand{\im}{{\rm im}\,}
\newcommand{\Pone}{$({P1})$}
\newcommand{\Ptwo}{$({P2})$}
\newcommand{\wPtwo}{$({wP2})$}
\begin{document}
\title[Groups with Submultiplicative Spectrum]{Finite Groups with Submultiplicative
Spectra 
}
\author{L. Grunenfelder, T. Ko\v sir, M. Omladi\v c, and H. Radjavi}%
\address{L. Grunenfelder: Department of Mathematics,
The University of British Columbia, 1984 Mathematics Road
Vancouver B.C., Canada V6T 1Z2}
\email{luzius@math.ubc.ca}%
\address{T. Ko\v sir, and M. Omladi\v c~: Department of Mathematics,
University of Ljubljana, Jadranska 19, 1000 Ljubljana, Slovenia}
\email{tomaz.kosir@fmf.uni-lj.si,\ matjaz.omladic@fmf.uni-lj.si}%
\address{Radjavi~: Department of Pure Mathematics,
University of Waterloo, Waterloo, Ontario,
Canada N2L 3G1}%
\email{hradjavi@uwaterloo.ca}
\thanks{Research supported in part by the NSERC of Canada and
by the Research Agency of Slovenia.}%

\date{\today}%

\begin{abstract}
We study abstract finite groups with the property, called
property \shat, that all of their subrepresentations have
submultiplicative spectra. Such groups are necessarily nilpotent
and we focus on $p$-groups. $p$-groups with property \shat\  are
regular. Hence, a $2$-group has property \shat\  if and only if
it is commutative. For an odd prime $p$, all $p$-abelian groups
have property \shat, in particular all groups of exponent $p$
have it. We show that a $3$-group or a metabelian $p$-group
($p\ge 5$) has property \shat\ if and only if it is V-regular.
\end{abstract}
\maketitle

\section{Introduction}
In recent years a number of properties of matrix groups (and
semigroups) were studied (see e.g. \cite{Radj,RaRo}). We wish to
propose a program to explore implications these results on matrix
groups might have for the theory of abstract groups: Given a
property ($P$) of matrix groups we say that an abstract group $G$
has property ($\widehat{P}$) if all the finite-dimensional
(irreducible) subrepresentations of $G$ have property ($P$).  We
call a representation of a subgroup of $G$ a {\em
subrepresentation} of $G$. In this paper we commence our program
by studying the so-called property \s.

Assume that $F$ is an algebraically closed field of
characteristic zero. A matrix group
$\mathcal{G}\subseteq GL_n(F)$ (or matrix semigroup
$\mathcal{G}\subseteq M_n(F)$) has {\em submultiplicative
spectrum} or, in short, it has {\em property} \s\ if for each
pair $A,B\in\mathcal{G}$ every eigenvalue of the product $AB$ is
equal to a product of an eigenvalue of $A$ and an eigenvalue of
$B$. Such groups and semigroups were first studied by Lambrou,
Longstaff, and Radjavi \cite{LaLoRa}. (See also \cite{Kram,Kram2,Kram3,Omla}.) If $\mathcal{G}$ is an
irreducible group with property \s\ then it is  nilpotent and
essentially finite \cite[Thms. 3.3.4 and 3.3.5]{RaRo}, i.e.,
$\mathcal{G}\subseteq F^*\mathcal{G}_0$ for some finite nilpotent
group $\mathcal{G}_0$. Here $F^*$ is the group of invertible
elements in $F$. A group is nilpotent if and only if it is the
direct product of its Sylow $p$-groups. So it is not a
restriction to study only $p$-groups.

A finite group $G$ has {\em property} \shat\ if all its
irreducible subrepresentations have property \s. Such groups are
necessarily nilpotent and we focus on $p$-groups with the
property. We show that a $2$-group has property \shat\ if and
only if it is commutative. For an odd prime $p$ we show that all
$p$-abelian groups have property \shat. A $p$-group $G$ is called
$p$-{\em{abelian}} if $(xy)^p=x^py^p$ for all $x,y\in G$. In
particular, all groups of exponent $p$ have property \shat. We
characterize all the metabelian $p$-groups with property \shat.
We show that $3$-groups and metabelian $p$-groups ($p\ge 5$) have
property \shat\ if and only if they are V-regular. In our proofs
we use several results on abstract $p$-groups, in particular
those of Alperin \cite{Alpe1,Alpe2}, Mann
\cite{Mann,MannII,MannI}, and Weichsel
\cite{Weich1,Weich2,Weich3}.

Let us remark that it would be interesting to consider a weaker
property than \shat. Namely, the property, called \stilde, that
every irreducible representation of $G$ has property \s. 
We do not know whether the properties \shat\ and \stilde\ are
equivalent for a finite $p$-group.

\section{Preliminaries}

We assume throughout that $G$ is a finite group and that $F$ is
an algebraically closed field of characteristic $0$. We denote by
$|g|$ the order of an element $g\in G$. The exponent $e(G)$ of
the group $G$ is the least common multiple of these orders.
In particular, if $G$ is a $p$-group,
then $e(G)$ is the maximum of the orders of its elements.

If $\varrho:G\to GL_n(F)$ is a representation, then
$\mathcal{G}=\varrho(G)$ is a finite matrix group. For
$A\in\mathcal{G}$ we denote by $\sigma(A)$ the spectrum of $A$.
We say that $\mathcal{G}$ has {\em property} \s\ if
\begin{equation}\label{submul}
\sigma(AB)\subseteq\sigma(A)\sigma(B)=\{\lambda\mu;\
\lambda\in\sigma(A),\mu\in\sigma(B)\}.
\end{equation}
for all $A,B\in\mathcal{G}$. An (abstract) group $G$ has {\em property} \shat\ if all its
irreducible subrepresentations have property \s. Here we call a
representation of a subgroup $H$ of $G$ a subrepresentation of
$G$. Groups with property \shat\ are nilpotent \cite[Thm.
3.3.5]{RaRo}. It is well known that an irreducible representation of a
nilpotent group is equivalent to a monomial representation (see
e.g. \cite[Thm. 16, p. 66]{Serre}, \cite[Lemma 6, p.
207]{Supr} or \cite[Cor. 6.3.11]{Grov}). Since each representation of a finite group $G$ is
completely reducible it follows that a group $G$ has property
\shat\ if and only if all its subrepresentations have property
\s. We use this fact later in the proofs, e. g., in the proof of
Proposition \ref{p1.5}.

Recall that all irreducible representations of a finite abelian
group have degree $1$. Hence we have:

\begin{lem}\label{abel} Every finite abelian group has property
\shat.
\end{lem}

The subgroups in the lower central series of $G$ are denoted by
$G^{(i)}$, i.e. $G^{(0)}=G$, $G^{(1)}=G'=[G,G]$, and
$G^{(i)}=[G^{(i-1)},G]$ for $i\ge 2$. We write $c=c(G)$ for the
{\em class} of $G$, i.e., $c$ is the least integer such that
$G^{(c)}=1$.

The subgroup of $G$ generated by the subset $\{ x_1,x_2,\ldots
,x_k\}$ will be denoted by $\langle x_1,x_2,\ldots,x_k\rangle$.

If $H\leq G$ is a subgroup and $K\unlhd H$ is a normal subgroup
then the quotient $H/K$ is called a {\em section} of $G$.

\begin{lem}\label{l0.0}
If $G$ has property \shat\ then all its subgroups and quotients
have property \shat. Furthermore, all its sections have property
\shat.
\end{lem}

{\em Proof}. It is enough to show that property \shat\ is
inherited by 
quotients.
Suppose that $H\unlhd G$ and that $\varrho:G/H\to GL_n(F)$ is a
representation. Then $\hat{\varrho}:G\to GL_n(F)$ defined by
$\hat{\varrho}(g)=\varrho(gH)$ is a representation of $G$ called
the inflation representation (see \cite[p.2]{Coll}). By our
assumption $(G,\hat{\varrho})$ has property \s, and hence so does
$(G/H,\varrho)$. \hfill\qed\vskip 5pt

The following lemma is an easy consequence of a theorem of
Burnside \cite[Thm. 1.2.2]{RaRo}.

\begin{lem}\label{l0.1}
If $\mathcal{G}_j\subseteq GL_{n_j}(F)$, $j=1,2$, are two
irreducible matrix groups then
$\mathcal{G}_1\otimes\mathcal{G}_2\subseteq GL_{n_1n_2}(F)$ is
also
irreducible.
\end{lem}

\begin{lem}\label{l0.2}
If $\mathcal{G}_j\subseteq GL_{n_j}(F)$, $j=1,2$, are two matrix
groups with property \s\ then also
$\mathcal{G}_1\otimes\mathcal{G}_2\subseteq GL_{n_1n_2}(F)$ has
property \s.
\end{lem}

{\em Proof}. Observe that $\sigma(A\otimes
B)=\sigma(A)\sigma(B)$. \hfill\qed\vskip 5pt

\begin{cor}\label{c0.3}
If $G_1$ and $G_2$ are finite groups with property \shat\ then
also the direct product $G_1\times G_2$ has property \shat.
\end{cor}

Since each finite nilpotent group is a direct product of its
Sylow $p$-groups we can limit our attention to $p$-groups.

\begin{prop}\label{p3.1}
A finite group $G$ has property \shat\ if and only if for each
pair of elements $x,y\in G$ the subgroup  $\langle x,y\rangle$
has
property \shat.
\end{prop}

{\em Proof}. If $G$ has property \shat\ then by definition
every subgroup, in particular every two-generated subgroup, has
property \shat.

Conversely, assume that every two generated subgroup of $G$ has
property \shat. Let $\varrho:K\to GL_n(F)$ be an irreducible
representation of a subgroup $K\subseteq G$. Choose $x,y\in K$
and let $H=\langle x,y\rangle$. The restriction $\varrho:H\to
GL_n(F)$ is a representation of $H$. By assumption it has property
\s\  and thus
$\sigma(\varrho(x)\varrho(y))\subseteq\sigma(\varrho(x))\sigma(\varrho(y))$.
\hfill\qed\vskip 5pt

\section{The Power Structure of $p$-Groups with Property \shat}

Suppose that $G$ is a finite $p$-group of exponent $p^e$. Then
for $k=1,2,\ldots,e$
$$\triangle_k(G)=\{g\in G;\ g^{p^k}=1\}$$
is the set of all the elements of order dividing $p^k$, and
$$\nabla_k(G)=\{g\in G;g=h^{p^k}\ {\rm for\ some}\ h\in G\}$$
is the set of all $p^k$-th powers. We denote by $\Omega_k(G)$ the
subgroup generated by $\triangle_k(G)$ and by $\mho_k(G)$ the
subgroup generated by $\nabla_k(G)$.

A $p$-group $G$ has {\em property} \Pone\ if for all the sections
$H$ of $G$ and all $k$ we have
$$
\nabla_k(H)=\mho_k(H).
$$
A $p$-group $G$ has {\em weak property} \Ptwo\ -- denoted by
\wPtwo\ -- if
$$
\triangle_k(G)=\Omega_k(G)
$$
for $k= 1,2,\ldots,e$. A $p$-group $G$ has {\em property} \Ptwo\
if all sections of $G$ have property \wPtwo.

Properties \Pone\ and \Ptwo\ were introduced by Mann
\cite{MannI}. We refer to \cite{MannI,Xu} for further details.

\begin{prop}\label{p1.4}
If a matrix group $\mathcal{G}\subseteq GL_n(F)$ has property \s\
then it has property \wPtwo.
\end{prop}

{\em Proof}. The submultiplicativity condition
$\sigma(AB)\subseteq\sigma(A)\sigma(B)$ implies that the order
$|AB|$ divides $\max\{|A|,|B|\}$. Hence, if
$A,B\in\triangle_k(\mathcal{G})$ then also
$AB,A^{-1}\in\triangle_k(\mathcal{G})$. \hfill\qed\vskip 5pt

\begin{prop}\label{p1.5}
If a $p$-group $G$ has property \shat\ then it has property \Ptwo.
\end{prop}

{\em Proof}. Suppose that $K$ is a section of $G$. By Lemma
\ref{l0.0} it follows that $K$ has property \shat. Take a
faithful representation $\varrho:K\to GL_n(F)$, e.g. the regular
representation. It has property \s\ and by Proposition \ref{p1.4}
it has property \wPtwo. Hence, $G$ has property
\Ptwo.\hfill\qed\vskip 5pt

Next we prove the main result of this section and one of our main
results. We begin by recalling some definitions.

A $p$-group is called {\em regular} if for every pair $x,y\in G$
there is an element $z$ in the commutator group $<x,y>'$ such that
$$(xy)^p=x^py^pz^p.$$
Note that \cite[Satz III.10.8(g)]{HuppI} shows that the above
definition of a regular $p$-group is equivalent to the more
common
one \cite[p. 321]{HuppI}.

A regular group $G$ is called {\em V-regular} if any finite
direct product of copies of $G$ is regular. Not every regular
$p$-group is V-regular -- see Wielandt's example \cite[Satz
III.10.3(c)]{HuppI}. A $p$-group $G$ is V-regular if and only if
all the finite groups in the variety of $G$ are regular
\cite{Neum,Weich2}. For the definitions of a variety of groups
and a variety of a given group $G$ we refer to Hanna Neumann's
book \cite{Neum}. Further properties of regular $p$-groups can be
found in \cite{HuppI,Suzu}.

\begin{thm}\label{t6.01}
If a $p$-group has property \shat\ then it is regular. Moreover,
it is V-regular.
\end{thm}

{\em Proof}. Assume that $G$ is a $p$-group with property \shat\
and that the exponent of $G$ is equal to $p^e$. If $e=1$ then it
is regular by \cite[Satz III.10.2(d)]{HuppI}. Suppose that $e\ge
2$. Now $G$ and the cyclic group $C_{p^e}$ of order $p^e$ both
have property \shat. The direct product $G\times C_{p^e}$ has
property \shat\ by Corollary \ref{c0.3}, property \Ptwo\ by
Proposition \ref{p1.5}, and property \Pone\ by \cite[Cor.
4]{MannI}. Finally, Theorem 25 of \cite{MannI} implies that $G$
is regular. The group $G$ is V-regular since by Corollary
\ref{c0.3} the direct product of any finite number of copies of
$G$ has property \shat\ and thus it is regular.\hfill\qed\vskip
5pt

Properties of regular groups \cite[Satz III.10.3(a),(b)]{HuppI}
are used to prove the following corollaries.

\begin{cor}\label{pis2}
A $2$-group has property \shat\ if and only if it is abelian.
\end{cor}

{\em Proof}. If a $2$-group is regular then it is abelian
\cite[Satz III.10.3(a)]{HuppI}. The converse follows since, by
Lemma \ref{abel}, all abelian groups have property \shat.
\hfill\qed\vskip 5pt

\vskip 10pt

Let us point out that if $k=1$ then finite matrix 2-groups in
$GL_{2^k}(F)$ with property \s\ are always commutative
\cite{LaLoRa}; however, they need not be commutative if $k\ge 2$
\cite{Kram,Omla}. Moreover, finite irreducible matrix $2$-groups
in $GL_{2^k}(F)$ with property \s\ are constructed in \cite{Omla}
for $k=3$ and in \cite{Kram} for $k\ge 4$.

\begin{cor}\label{pis3}
If a $3$-group has property \shat\ then it is metabelian.
\end{cor}

{\em Proof}. This follows from a result of Alperin \cite[Thm.
1]{Alpe1}.\hfill\qed\vskip 5pt

\begin{prop}\label{p6.01}
If $G$ has property \shat\ then any finite group in the variety
of
$G$ has property \shat.
\end{prop}

{\em Proof}. By \cite[Cor. 32.32]{Neum} a finite group $H$ in the
variety of $G$ is a section of a finite direct product of copies
of $G$. Corollary \ref{c0.3} implies that the direct product has
property \shat\ and Lemma \ref{l0.0} implies that $H$ has
property \shat\ as well.\hfill\qed\vskip 5pt\vskip 10pt

In \S{6} we show that the converse of Theorem \ref{t6.01} is true
for metabelian $p$-groups. We do not know the answer to the
general question: Does a finite V-regular $p$-group have property
\shat? We know that the direct product of finitely many groups
with property \shat\ again has property \shat. If the direct
product of two V-regular groups is not V-regular then the answer
to the above question is negative. The question if the direct
product of two V-regular groups is V-regular was studied by
Groves \cite{Groves2}.

\section{Matrix Groups in $GL_p(F)$ with Property \s}

In this section we consider an irreducible matrix $p$-group
$\mathcal{G}$ in $GL_p(F)$. We assume hereafter that $p$ is an odd
prime. The main result of the section is the following: if
$\mathcal{G}$ has property \s\ then the class $c(\mathcal{G})$ is
at most $p-1$.

Assume first that $\mathcal{G}\subseteq GL_p(F)$ is an irreducible
$p$-group of exponent $e(G)=p$. Then we may assume without loss that
$\mathcal{G}$ is monomial. Each element of $\mathcal{G}$ is either diagonal or
of the form $DP^k$, where $D$ is diagonal, $k\in\{1,2,\ldots,p-1\}$ and
$$P=\left[\begin{array}{cccccc}
  0 & 0 & 0 & \cdots & 0 & 1 \\
  1 & 0 & 0 & \cdots & 0 & 0 \\
  0 & 1 & 0 & \cdots & 0 & 0 \\
  \vdots & \vdots & \ddots & \ddots & \vdots & \vdots \\
  0 & 0 & 0 & \cdots & 0 & 0 \\
  0 & 0 & 0 & \cdots & 1 & 0
\end{array}\right].$$
Since $e(G)=p$ and $(DP^k)^p=(\det D) I$ it follows that $\det D=1$.
Note that an element of the form $DP$, $D$ diagonal of determinant $1$, is diagonally similar
to $P$. Therefore we may assume without loss that $P\in\mathcal{G}$. We denote by $ \mathcal{D}$ the subgroup of
all the diagonal elements in $ \mathcal{G}$. A simple matrix computation shows  that
$G'\subseteq\mathcal{D}$. Let $\omega$ be a primitive $p$-th root
of $1$ and $\Gamma_1$ the set of all the $p$-th roots of $1$.
Further we denote by
$\mathbb{Z}_p=\mathbb{Z}/p\mathbb{\mathbb{Z}}$ the finite field
with $p$ elements. We define a map $\chi:\mathcal{D}\to
\mathbb{Z}_p^p$ by
$$\chi\left[\begin{array}{cccccc}
  \omega^{k_1} & 0 & 0 & \cdots & 0  \\
  0 & \omega^{k_2} & 0 & \cdots & 0  \\
  0 & 0 & \omega^{k_3} & \cdots & 0  \\
  \vdots & \vdots & \ddots & \vdots & \vdots \\
  0 & 0 & 0 & \cdots & \omega^{k_p}
\end{array}\right]=(k_1,k_2,\ldots,k_p).$$

\begin{lem}\label{l1}
$\chi$ is a homomorphism of abelian groups and its image $\im
\chi$ is invariant under the cyclic permutation
$\pi:\mathbb{Z}_p^p\to\mathbb{Z}_p^p$ given by
$$\pi(k_1,k_2,\ldots,k_p)=(k_2,k_3,\ldots,k_p,k_1).$$
\end{lem}

{\em Proof}. It is an easy observation that $\chi$ is a
homomorphism and that its image is a vector subspace. Since
$$P^{-1}\left[\begin{array}{cccccc}
  \omega^{k_1} & 0 & 0 & \cdots & 0  \\
  0 & \omega^{k_2} & 0 & \cdots & 0  \\
  0 & 0 & \omega^{k_3} & \cdots & 0  \\
  \vdots & \vdots & \ddots & \vdots & \vdots \\
  0 & 0 & 0 & \cdots & \omega^{k_p}
\end{array}\right]P=\left[\begin{array}{cccccc}
  \omega^{k_2} & 0  & \cdots & 0 & 0 \\
  0 & \omega^{k_3}  & \cdots & 0 & 0  \\
  \vdots & \vdots & \ddots & \vdots & \vdots \\
  0 & 0 & \cdots & \omega^{k_p} & 0 \\
  0 & 0 & \cdots & 0 & \omega^{k_1}
\end{array}\right]$$
it follows that $\pi(\im\chi)\subseteq\im\chi$.\hfill\qed\vskip
5pt

\begin{lem}\label{l2}
There are exactly $p+1$ subspaces in $\mathbb{Z}_p^p$ invariant
under $\pi$, one in each dimension $j=0,1,\ldots,p$. They are
$\im(I-\pi)^{p-j}$ for $j=0,1,\ldots,p-1$, and $\mathbb{Z}_p^p$.
Also, $(I-\pi)^p=0$.
\end{lem}

{\em Proof}. The linear maps $\pi$ and $I-\pi$ have the same
invariant subspaces. Since $\pi^p=I$ it follows that
$(I-\pi)^p=0$. The matrix
$$I-\pi=\left[\begin{array}{ccccc}
  1 & 0 & 0 & \cdots & -1 \\
  -1 & 1 & 0 & \cdots & 0 \\
  0 & -1 & 1 & \cdots & 0 \\
  \vdots & \vdots & \ddots & \ddots & \vdots \\
  0 & 0 & 0 & \cdots & 1
\end{array}\right]$$
has rank equal to $p-1$. Hence
$$
\mathbb{Z}_p^p\supset\im(I-\pi)\supset\im(I-\pi)^2\supset\cdots
\supset\im(I-\pi)^{p-1}\supset0$$ is the chain of all the
distinct
invariant subspaces of $I-\pi$.\hfill\qed\vskip 5pt

\begin{lem}\label{l3}
For $j\ge 1$ we have
$\chi\left(\mathcal{G}^{(j)}\right)\subseteq\im (I-\pi)^j.$
\end{lem}

{\em Proof}. Since $\mathcal{G}$ is monomial and $P\in\mathcal{G}$ it follows that each element of $ \mathcal{G}$ can be
written in the form
\begin{equation}\label{PD}
P^lD_1=D_2P^l
\end{equation}
for some $l\in\{0,1,\ldots,p-1\}$  and $D_1,D_2\in\mathcal{D}$.

We prove the lemma by induction on $j$. Assume $j=1$. It is an
easy consequence of the form (\ref{PD}) that elements of
$\mathcal{G}^{(1)}$ are products of elements of the form
$DP^lD^{-1}P^{-l}$ for some $l\in\{1,\ldots,p-1\}$  and
$D\in\mathcal{D}$. If
$$D=\left[\begin{array}{cccccc}
  \omega^{k_1} & 0 & 0 & \cdots & 0  \\
  0 & \omega^{k_2} & 0 & \cdots & 0  \\
  0 & 0 & \omega^{k_3} & \cdots & 0  \\
  \vdots & \vdots & \ddots & \vdots & \vdots \\
  0 & 0 & 0 & \cdots & \omega^{k_p}
\end{array}\right]$$
then
$$DP^lD^{-1}P^{-l}=\left[\begin{array}{cccccc}
  \omega^{k_1-k_{l+1}} & 0 & 0 & \cdots & 0  \\
  0 & \omega^{k_2-k_{l+2}} & 0 & \cdots & 0  \\
  0 & 0 & \omega^{k_{l+3}} & \cdots & 0  \\
  \vdots & \vdots & \ddots & \vdots & \vdots \\
  0 & 0 & 0 & \cdots & \omega^{k_p-k_l}
\end{array}\right],$$
where the index $s$ of $k_s$ is computed modulo $p$. It follows
that $\chi([D,P^l])\subseteq\im(I-\pi^l)$. Since
$I-\pi^l=(I-\pi)(I+\pi+\pi^2+\cdots+\pi^{l-1})$ we see that
$\im(I-\pi^l)\subseteq\im(I-\pi)$ for $l=1,2,\ldots,p-1$.
Therefore, $\chi\left(\mathcal{G}^{(1)}\right)
\subseteq\im(I-\pi)$.

Assume now that $D\in\mathcal{G}^{(j-1)}$. The induction
hypothesis is that $\chi(D)\in\im(I-\pi)^{j-1}$. An easy matrix
computation shows that each element of $\mathcal{G}^{(j)}$ is a product of
elements of the form $DP^lD^{-1}P^{-l}$ for some
$l\in\{1,\ldots,p-1\}$ and $D\in\mathcal{G}^{(j-1)}$. Then we
prove, in a way similar to the case $j=1$, that
$\chi(DP^lD^{-1}P^{-l})\in\im(I-\pi)^j$ and thus
$\chi\left(\mathcal{G}^{(j)}\right)\subseteq\im(I-\pi)^j$.\hfill\qed\vskip
5pt

\begin{cor}\label{c1}
If $\mathcal{G}\subseteq GL_p(F)$ is an irreducible $p$-group of
exponent $p$ then its class is at most $p-1$.
\end{cor}

{\em Proof}. Lemma \ref{l2} and Lemma \ref{l3} with $l=p$ imply
that $\chi(\mathcal{G}^{(p)})\subseteq\im(I-\pi)^p=0$. Therefore,
$\mathcal{G}^{(p)}=1$.\hfill\qed\vskip 5pt

\begin{prop}\label{p4}
If $\mathcal{G}\subseteq SL_p(F)$ is an irreducible $p$-group with
property \wPtwo\ then the exponent of $\mathcal{G}$ is equal to
$p$.
\end{prop}

{\em Proof}. We denote by $\mathcal{D}$ the subgroup of all the
diagonal elements of $\mathcal{G}$. Assume that
$A\in\mathcal{G}\backslash\mathcal{D}$. Then $A^p=I$ since $\det
A=1$. We may assume that
$$P=\left[\begin{array}{cccccc}
  0 & 0 & 0 & \cdots & 0 & 1 \\
  1 & 0 & 0 & \cdots & 0 & 0 \\
  0 & 1 & 0 & \cdots & 0 & 0 \\
  \vdots & \vdots & \ddots & \ddots & \vdots & \vdots \\
  0 & 0 & 0 & \cdots & 0 & 0 \\
  0 & 0 & 0 & \cdots & 1 & 0
\end{array}\right]$$
is in $\mathcal{G}$. Suppose now that $D\in\mathcal{D}$. Since
$\det D=1$ it follows that $(DP)^p=1$. Hence
$P,DP\in\Omega_1(\mathcal{G})$. Since $\mathcal{G}$ has property
\wPtwo\ it follows that $D\in\Omega_1(\mathcal{G})$ and hence
$D^p=I$.\hfill\qed\vskip 5pt

\begin{thm}\label{t2}
If $\mathcal{G}\subseteq GL_p(F)$ is an irreducible $p$-group with
property \s\ then its class is at most $p-1$.
\end{thm}

{\em Proof}. Assume that the exponent of $\mathcal{G}$ is equal
to $p^e$. Let $\theta$ be a primitive $p^{e+1}$-th root of $1$.
We enlarge $\mathcal{G}$ to ${\mathcal{H}}$ by multiplying all
the elements by $\theta^j$, $j=1,2,\ldots,p^{e+1}$. Note that
this only enlarges the center, all other quotients of two
consecutive elements of the upper central series of $\mathcal{G}$
and ${\mathcal{H}}$ are equal. Hence the classes of both groups
are equal. Next we consider the subgroup
$\mathcal{K}=\left\{A\in\mathcal{H};\ \det A=1\right\}$. Since
the exponent of $\mathcal{G}$ is equal to $p^e$ it follows that
for each $A\in\mathcal{G}$ there is an integer $k(A)$ such that
$\theta^{k(A)}A\in\mathcal{K}$. The elements of $\mathcal{G}'$
are products of commutators $[A,B]$. Note that each commutator
$[A,B]$ has determinant equal to $1$. Since
$[A,B]=[\theta^{k(A)}A,\theta^{k(B)}B]$ it follows that
$\mathcal{G}^{(j)}=\mathcal{K}^{(j)}$ for $j=1,2,\ldots$, and
hence the classes $c(\mathcal{G})$ and $c(\mathcal{K})$ are
equal. By Proposition \ref{p1.4} property \wPtwo\ follows from
property \s. Next, Proposition \ref{p4} implies that the exponent
of $\mathcal{K}$ is equal to $p$ and Corollary \ref{c1} implies
that $c(\mathcal{K})\le p-1$. \hfill\qed\vskip 5pt

\section{$p$-Abelian Groups Have Property \shat}

\begin{thm}\label{expp}
If the exponent of $G$ is equal to $p$ then $G$ has property
\shat.
\end{thm}

{\em Proof}. It suffices to show that each finite irreducible
matrix group $\mathcal{G}\subseteq GL_{p^k}(F)$, $k\ge 0$, of
exponent $p$ has property \s. Assume that $\mathcal{G}$ is
monomial and denote by $\mathcal{D}$ the subgroup of all the
diagonal matrices. Since the exponent of $\mathcal{G}$ is equal
to $p$ it follows that $\sigma(D)\subseteq\Gamma_1$ for all
$D\in\mathcal{D}$. Each element of $\mathcal{G}$ is of the form
$DP$ for some $D\in\mathcal{D}$ and a permutation matrix $P$ of
order dividing $p$.

We choose two elements $A_1=D_1P_1$ and $A_2=D_2P_2$ in $\mathcal{G}$.
Here $D_1,D_2\in\mathcal{D}$ and $P_1,P_2$ are permutation
matrices. Observe that our assumptions imply that if $P_i\neq I$
then $\sigma(A_i)=\Gamma_1$.

To show submultiplicativity of spectra we treat three cases:
\begin{itemize}
  \item If $P_1=P_2=I$ then the submultiplicativity is obvious.
  \item If $P_1P_2\neq I$ then one of $P_1,P_2$
  is not equal to $I$. We assume $P_1\neq I$, the case $P_1=I$ and $P_2\neq
  I$ is done in a similar way. Then
$$\sigma(A_1A_2)=\Gamma_1=\Gamma_1\sigma(A_2)=\sigma(A_1)\sigma(A_2).$$
  \item If $P_1P_2=I$, but neither of $P_1,P_2$ is equal to
  $I$, then
$$\sigma(A_1A_2)\subseteq\Gamma_1=\Gamma_1\Gamma_1=\sigma(A_1)\sigma(A_2).$$
\end{itemize}
\hfill\qed\vskip 5pt

A finite $p$-group $G$ is called $p$-abelian if $(xy)^p=x^py^p$
for all $x,y\in G$. Now we recall a characterization of such
groups \cite{Alpe2,Weich3}: A finite group is $p$-abelian if and
only if it is a section of a direct product of an abelian
$p$-group and a group of exponent $p$. By Lemma \ref{abel}
abelian groups have property \shat. So we have the following
consequence of Corollary \ref{c0.3} and Theorem \ref{expp}.

\begin{cor}\label{pabelian}
A finite $p$-abelian group has property \shat.
\end{cor}

The following result is of interest on its own, and it will be
used later as the first step of a proof by induction.

\begin{cor}\label{det1A}
Suppose $\mathcal{G}\subseteq SL_p(F)$ is an irreducible $p$-group. Then the
following are equivalent:
\begin{enumerate}
  \item $\mathcal{G}$ has property \s,
  \item $\mathcal{G}$ has property \wPtwo,
  \item $e(\mathcal{G})=p$.
\end{enumerate}
\end{cor}

{\em Proof}. The implication $(1)\Rightarrow (2)$ follows by
Proposition \ref{p1.4}, the implication $(2)\Rightarrow (3)$ by
Proposition \ref{p4}, and $(3)\Rightarrow (1)$ by Theorem
\ref{expp}. \hfill\qed\vskip 5pt

\section{Metabelian Groups with Property \shat}

In this section we assume that $G$ is a {\em metabelian}
$p$-group, i.e. we assume that $G'$ is abelian. Our result
extends a result of Weichsel \cite{Weich2} that characterizes
metabelian V-regular $p$-groups. We show that a finite metabelian
group has property \shat\ if and only if it is V-regular.

Let us introduce some notation. We write
$$P_k=
\left[\begin{array}{ccccc}
  0 & 0 & \cdots & 0 & 1 \\
  1 & 0 & \cdots & 0 & 0 \\
  0 & 1 & \cdots & 0 & 0 \\
  \vdots & \vdots & \ddots & \ddots & \vdots \\
  0 & 0 & \cdots & 1 & 0
\end{array}\right]$$
for the cyclic matrix of order $p^k$ in $GL_{p^k}(F)$. If $D\in
GL_{p^k}(F)$ is a diagonal matrix of order $p^l$ for some $l$
then an element of the form $DP_k$ is called a {\em big cycle}.

\begin{lem}\label{det1}
Suppose that a monomial $p$-group $\mathcal{G}\subseteq SL_{p^k}(F)$ is
generated by a big cycle and a diagonal matrix. If $\mathcal{G}$
is irreducible with property \wPtwo\ then the exponent of $\mathcal{G}$
is equal to $p^k$.
\end{lem}

{\em Proof}. Assume that the generators of $\mathcal{G}$ are a
big cycle $DP_k$ and a diagonal matrix $B$. Here $D$ is a
diagonal matrix, too. Since $\det(DP_k)=\det D=1$ it follows that
$DP_k$ is similar  to $P_k$ using a diagonal similarity. Without loss we may assume that
$D=I$, i.e., that $P_k\in \mathcal{G}$.
We denote by $\mathcal{D}$ the subgroup of
all the diagonal matrices in $\mathcal{G}$. Then each element of
$\mathcal{G}$ is of the form $EP_k^j$ for some $E\in\mathcal{D}$
and some integer $j$.

We prove the lemma by induction on $k$. The
case $k=1$ was proved in Proposition \ref{p4}. Assume now that
$k\ge 2$. Suppose that the subgroup
$\mathcal{H}\subseteq\mathcal{G}$ consists of all the elements of
the form $EP_k^j$, where $E\in\mathcal{D}$ and $j$ is a multiple
of $p$.  We may assume that, up to a permutational similarity,
the elements of $\mathcal{H}$ are all of the form
\begin{equation}\label{Aform}
A=\left[\begin{array}{cccc}
  A_1 & 0 & \cdots & 0 \\
  0 & A_2 & \cdots & 0 \\
  \vdots & \vdots & \ddots & \vdots \\
  0 & 0 & \cdots & A_{p}
\end{array}\right],
\end{equation}
where $A_j\in GL_{p^{k-1}}(F)$, $j=1,2,\ldots,p$. We denote by $\mathcal{H}_1$ the
subgroup in $GL_{p^{k-1}}(F)$ generated
by all the blocks $A_1$ of all elements $A\in\mathcal{H}$. Observe that it is
irreducible. Let
$$
\widetilde{\mathcal{H}}_1=\left\{\theta B;\ det (\theta
B)=1,\theta\in F,B\in\mathcal{H}_1\right\}.
$$
Then the group $\widetilde{\mathcal{H}}_1$ is an irreducible
$p$-group in $GL_{p^{k-1}}(F)$ such that $\det B=1$ for all
$B\in\widetilde{\mathcal{H}}_1$. By the inductive hypothesis the
exponent $e(\widetilde{\mathcal{H}}_1)$ is equal to $p^{k-1}$. Choose
an element $C\in\mathcal{G}\backslash\mathcal{H}$. Without loss
we may assume that
$$
C=\left[\begin{array}{ccccc}
  0 & 0 & \cdots & 0 & U \\
  I & 0 & \cdots & 0 & 0 \\
  0 & I & \cdots & 0 & 0 \\
  \vdots & \vdots & \ddots & \ddots & \vdots \\
  0 & 0 & \cdots & I & 0
\end{array}\right],
$$
where $U\in GL_{p^{k-1}}(F)$. Observe that $\det U=1$
since $\det C=1$ and that $C^{p}\in\mathcal{H}$. Hence $U\in\widetilde{\mathcal{H}}_1$.
Then $|U|$ divides $p^{k-1}$ and $|C|$ divides $p^k$. Next assume
that $A\in\mathcal{H}$ is of form (\ref{Aform}). Since $\det A=1$
it follows that $\prod_{j=1}^{p}\det A_j=1$. In the same way as
we did for $C$ we prove that $|AC|$ divides $p^k$. Since
$\mathcal{G}$ has property \wPtwo\ it follows that $|A|$ also
divides $p^k$.
This shows that $e(\mathcal{G})$ divides $p^k$. Since $|P_k|=p^k$ it follows that
$e(\mathcal{G})=p^k$.\hfill\qed\vskip 5pt

We denote by $\Gamma_k$ the set of all $p^k$-th roots of $1$.
If $\eta\in\Gamma_k$ is a scalar and $i$ a positive integer then
$$
D_k(i,\eta)= \left[\begin{array}{ccccc}
  1 & 0 & 0 & \cdots & 0 \\
  0 & \eta^{{1\choose i}} & 0 & \cdots & 0 \\
  0 & 0 & \eta^{{2\choose i}} & \cdots & 0 \\
  \vdots & \vdots & \vdots & \ddots & \vdots \\
  0 & 0 & 0 & \cdots & \eta^{{p^k-1\choose
  i}}\end{array}\right]
$$
is a diagonal matrix in $GL_{p^k}(F)$. Here we assume that
${j\choose i}=0$ if $j<i$.

\begin{lem}\label{centers}
Assume that a $2$-generated irreducible monomial $p$-group
$\mathcal{G}\subseteq SL_{p^k}(F)$ has
class $c\le p-1$. Suppose further that one of the generators is a big
cycle and the other is a diagonal matrix. Then $\mathcal{G}$ has property \wPtwo, its exponent is equal to $p^k$,
and each element of the $(c-i)$-th subgroup $\mathcal{G}^{(c-i)}$ in the lower central series of $\mathcal{G}$,
$i=1,2,\ldots,c-1,$ is a product of elements of the form
\begin{equation}\label{generators}
\alpha_0 I,\ and\ \alpha_j D_k(j,\eta_j),\ j=1,2,\ldots,i-1,
\end{equation}
for some $\alpha_0\in F$, $\alpha_j\in\Gamma_k$, $\eta_j\in\Gamma_{k}$, $j=1,2,\ldots,i-1$.
\end{lem}

{\em Proof}. Since $c\le p-1$ it follows that $\mathcal{G}$ is a regular group \cite[p. 322]{HuppI}, and hence it has
properties \Ptwo\ and \wPtwo\ \cite{MannI}. By Lemma \ref{det1} its exponent is equal to $p^k$.
The irreducibility of $\mathcal{G}$ implies that its center
consists of scalar matrices, which have order dividing $p^k$. Assume that $B$ is the diagonal generator and the
other generator is $C=DP_k$, where $D$ is a diagonal matrix. Since $\det A=1$ for all $A\in \mathcal{G}$
it follows that $\det C=\det D =1$. So $C$ is similar to $P_k$ using a diagonal similarity. Without loss
we may assume that $C=P_k\in\mathcal{G}$.
We denote by $\mathcal{D}$ the subgroup of all diagonal matrices in $\mathcal{G}$. Since $\mathcal{G}$ is monomial and $P_k\in\mathcal{G}$ it follows that each element of $ \mathcal{G}$ can be
written in the form
\begin{equation}\label{PDdav}
P_k^lD_1=D_2P_k^l
\end{equation}
for some $l\in\{0,1,\ldots,p^k-1\}$  and $D_1,D_2\in\mathcal{D}$.
It is an easy consequence of the form (\ref{PDdav}) that elements of
$\mathcal{G}^{(1)}$ are products of elements of the form
$DP_k^lD^{-1}P^{-l}_k$ for some $l\in\{1,\ldots,p^k-1\}$  and
$D\in\mathcal{D}$. In particular, it follows that
$\mathcal{G}^{(1)}\subset\mathcal{D}$. Observe that $\mathcal{G}^{(c-1)}$ is a nontrivial subgroup of the center
$Z(\mathcal{G})$ of $\mathcal{G}$.

Let $\omega=e^{{\frac{2\pi i}{2^k}}}$ be a primitive $p^k$-th root
of $1$ and thus $\Gamma_k=\{\omega^j,\ j=0,1,\ldots,p^k-1\}$.
Further we denote by
$\mathbb{Z}_{p^k}=\mathbb{Z}/{p^k}{\mathbb{Z}}$ the finite quotient ring of $\mathbb{Z}$ by the principal ideal generated by $p^k$. We define a map $\chi:\mathcal{D}\to
\mathbb{Z}_{p^k}^{p^k}$ by
$$\chi\left[\begin{array}{cccccc}
  \omega^{l_1} & 0 & 0 & \cdots & 0  \\
  0 & \omega^{l_2} & 0 & \cdots & 0  \\
  0 & 0 & \omega^{l_3} & \cdots & 0  \\
  \vdots & \vdots & \ddots & \vdots & \vdots \\
  0 & 0 & 0 & \cdots & \omega^{l_{p^k}}
\end{array}\right]=(l_1,l_2,\ldots,l_{p^k}).$$
A proof similar to the proof of Lemma \ref{l1} shows that
$\chi$ is a homomorphism of abelian groups.
We define the cyclic permutation
$\pi:\mathbb{Z}_{p^k}^{p^k}\to\mathbb{Z}_{p^k}^{p^k}$ by
$$\pi(l_1,l_2,\ldots,l_{p^k})=(l_2,l_3,\ldots,l_{p^k},l_1).$$

If a matrix $C_2$ in $\mathcal{G}^{(c-2)}$ is such that $[P_k,C_2]\neq I$ then
$$[P_k,C_2]\in\mathcal{G}^{(c-1)}\subset Z(\mathcal{G})$$
and so
\begin{equation}\label{C2}
[P_k,C_2]=\omega^t I\
\end{equation}
for some $t$ such that $1\le t\le p^k-1$. If we write
$$C_2=
\left[\begin{array}{ccccc}
  \omega^{l_1} & 0 & 0 & \cdots & 0 \\
  0 & \omega^{l_2} & 0 & \cdots & 0 \\
  0 & 0 & \omega^{l_3} & \cdots & 0 \\
  \vdots & \vdots & \vdots & \ddots & \vdots \\
  0 & 0 & 0 & \cdots & \omega^{l_{p^k}}\end{array}\right]
$$
then  (\ref{C2}) implies that
\begin{equation}\label{diffeq}
l_{j+1}-l_j=t,\ j=1,2,3,\ldots,p^k-1
\end{equation} and
\begin{equation}\label{bc}
l_1-l_{p^k}=t.
\end{equation}
This can be viewed as a simple linear difference equation (\ref{diffeq}) for an infinite
sequence $\{l_j\}_{j=1}^{\infty}$. Its solution is of the form
\begin{equation}\label{solution1}
l_{j}=t_1 j + t_0 = t_1 {j\choose 1} + t_0 {j\choose 0}, \ j=1,2,3,\ldots,p^k-1,
\end{equation}
for some $t_1,t_0\in\mathbb{Z}$.
Observe that condition (\ref{bc}) is satisfied modulo $p^k$.
Using (\ref{solution1}) we obtain that $C_2=\alpha D_k(1,\eta)$ for some
scalars $\alpha,\eta\in\Gamma_k$. Note that the expression for $l_j$ in (\ref{solution1}) is linear in $j$.

We prove the structure result for elements in $\mathcal{G}^{(c-i)}$ by induction on $i$. Our inductive assumption is that
for each $C_i\in\mathcal{G}^{(c-i)}$ the elements $l_j$ of the sequence $\chi(C_i)=(l_j)_{j=1}^{p^k}$ are given as a linear combination of binomial expressions ${j\choose u}$, $u=0,1,\ldots,i-1$, with integer coefficients. The case
$i=2$ was proved above.

Now we take an element $C_{i+1}\in \mathcal{G}^{(c-i-1)}$ and denote by
$(l_j)_{j=1}^{p^k}$ the image $\chi(C_{i+1})$. Then $[C_{i+1},P_k]$ is in $\in\mathcal{G}^{(c-i)}$ and we have
\begin{equation}\label{C_i+1}
\chi([C_{i+1},P_k])=(I-\pi)\chi(C_{i+1}).
\end{equation}
The inductive assumption implies that the elements of (\ref{C_i+1}) are given as a linear combination of binomial expressions ${j\choose u}$, $u=0,1,\ldots,i-1$, with integer coefficients. As before, we can view the components of (\ref{C_i+1}) as a simple linear difference equation.
Its solution, i.e. the elements of $\chi(C_{i+1})$ are then given by
\begin{equation}\label{solution2}
l_{j}=\sum_{u=0}^i s_u {j\choose u}, \ j=1,2,3,\ldots,p^k.
\end{equation}
Since any polynomial which has integer values if the argument is an integer can be written as a $\mathbb{Z}$-linear combination in the binomial basis ${j\choose u}$ (confer \cite[p. 2]{CaCh}) it follows that the coefficients $s_u$ in (\ref{solution2}) are integers.
Then
\begin{equation}\label{addp_k}
l_{p^m+j}-l_j=\sum_{u=0}^i s_u \left( {p^m+j\choose u} - {j\choose u}\right)
\end{equation}
for $m=1,2,\ldots$ Since $0\le u\le i<p$, it is clear that
$${p^m+j\choose u} - {j\choose u}$$
is divisible by $p^m$ and hence $p^m$ divides $l_{p^m+j}-l_j$. In particular,
$$l_{p^k}-l_1=\left(l_{p^k}-l_{p^k+1}\right)+\left(l_{p^k+1}-l_1\right)\equiv\left(l_{p^k}-l_{p^k+1}\right)\quad {\pmod{p^k}}.$$

This implies that $l_{p^k+1}-l_1$ is divisible by $p^k$. In particular, this implies that also the equation given by the $p^k$-th component of (\ref{C_i+1}) is satisfied modulo $p^k$.
Finally, since $l_j$ can be written in the basis given by the binomial expressions ${j\choose u}$ it follows that $C_{i+1}$
is a product of
elements of the forms
$$\alpha_0 I,\ \alpha_jD_k(j,\eta_j),\ j=1,2,\ldots,i,$$
where $\alpha_0,\ \alpha_j, \eta_j\in\Gamma_{k}$, $j=1,2,\ldots,i$.
\hfill\qed\vskip 5pt

\begin{lem}\label{modp}
The matrices $D_k(i,\eta)$, for $i=1,2,\ldots,p-2$ and
$\eta\in\Gamma_k$, have the following properties:
\begin{enumerate}
  \item $\det D_k(i,\eta)=1$,
  \item $D_k(i,\eta)$ is permutationally similar to
$$
\widetilde{D}_k(i,\eta)=\left[\begin{array}{cccc}
  \alpha_1 E_1 & 0 & \cdots & 0 \\
  0 & \alpha_2 E_2 & \cdots & 0 \\
  \vdots & \vdots & \ddots & \vdots \\
  0 & 0 & \cdots & \alpha_p E_{p}
\end{array}\right],
$$
where the matrices $E_1,E_2,\ldots,E_p$ are diagonal with the
determinant equal to $1$ and each is a product of elements of the form
$D_{k-1}(i,\theta)$ for some scalars $\theta\in\Gamma_{k-1}$. The similarity
between $D_k(i,\eta)$ and $\widetilde{D}_k(i,\eta)$ is induced by
the reordering of the standard basis $(e_1,e_2,\ldots,e_{p^k})$
to $(e_1$, $e_{p^{k-1}+1}$, $\ldots$~, $e_{(p-1)p^{k-1}+1}$,
$e_2$, $e_{p^{k-1}+2}$, $\ldots$~, $e_{(p-1)p^{k-1}+2}$,
$\ldots$~, $e_p$, $e_{p^{k-1}+p}$, $\ldots$~,$e_{p^k})$.
\end{enumerate}
\end{lem}

{\em Proof}. Property (1) follows from the identity
$$
\sum_{j=1}^n {j\choose i}={{n+1}\choose {i+1}},
$$
which holds for all positive integers $i$ and $n$ and can be verified by a counting argument.

Property (2) follows from the fact that the elements of the sequence $\chi(D_k(i,\eta))$ are given by an expression of the form (\ref{solution2}), which satisfies relation (\ref{addp_k}). Taking $m=1$ we see that $E_j$ are products of elements of the form $D_{k-1}(i,\theta)$ for some scalars $\theta\in\Gamma_{k-1}$.
\hfill\qed\vskip 5pt

\begin{prop}\label{bigcycle}
Suppose that $\mathcal{G}\subseteq SL_{p^k}(F)$ is an irreducible monomial
$p$-group that is generated by a big cycle and a
diagonal matrix. If $\mathcal{G}$ has class at most $p-1$ then
it has property \s.
\end{prop}

{\em Proof}. Assume that $DP_k$ is the big cycle generator, where
$D$ is a diagonal matrix. Since $1=\det DP_k=\det D$ it follows
that $DP_k$ is similar, in fact by a diagonal similarity, to
$P_k$. Thus, we may further assume that $P_k\in\mathcal{G}$. We
denote by $\mathcal{D}$ the subgroup of all the diagonal matrices
in $\mathcal{G}$. Each element of $\mathcal{G}$ is of the
form $DP^j_k$ for some integer $j$ and matrix $D\in\mathcal{D}$.

Observe that Lemma  \ref{centers} implies that $\mathcal{G}$ has
exponent equal to $p^k$ and it has property \wPtwo.
We prove the proposition by induction on $k$. For
$k=1$ the claim follows by Corollary \ref{det1A}. Assume that our
claim is true for the subgroups of $SL_{p^l}(F)$ with $l<k$.
Choose two elements $A_1=D_1P_k^{j_1}$ and $A_2=D_2P_k^{j_2}$ in
$\mathcal{G}$. We consider several cases:
\begin{itemize}
  \item If $j_1=j_2=0$ then the submultiplicativity is obvious.
  \item If $j_1+j_2$ is not divisible by $p$ then one of $j_1,j_2$
  is not divisible by $p$. We assume that $j_1$ is not divisible
  by $p$. (The case $j_1$ is divisible by $p$ and $j_2$
  is not divisible by $p$ is done in a similar way.) Then
$$\sigma(A_1A_2)=\Gamma_k=\Gamma_k\sigma(A_2)=\sigma(A_1)\sigma(A_2).$$
  \item If $j_1+j_2$ is $0$ modulo $p^k$, but neither of $j_1,j_2$ is $0$ or divisible by $p$, then
$$\sigma(A_1A_2)\subseteq\Gamma_k=\Gamma_k\Gamma_k=\sigma(A_1)\sigma(A_2).$$
\end{itemize}

It remains to consider the case when both $j_1$ and $j_2$ are
divisible by $p$. By Lemma \ref{centers} it follows that each element of
the subgroup $\mathcal{G}^{(c-i)}$
$i=1,2,\ldots,c-1,$ is equal to a product of elements of the following possible forms:
\begin{equation*}
\beta_0I,\ \beta_j D_k(j,\eta_j),\ j=1,2,\ldots,i-1,
\end{equation*}
where $\beta_j\in\Gamma_k$ and $\eta_j\in\Gamma_{k}$.
By Lemma \ref{modp} each element of the form $D_k(j,\eta_j)$ for $j\le p-2$, is
permutationally similar to a matrix of the form
$$
\left[\begin{array}{cccc}
  \alpha_1 E_1 & 0 & \cdots & 0 \\
  0 & \alpha_2 E_2 & \cdots & 0 \\
  \vdots & \vdots & \ddots & \vdots \\
  0 & 0 & \cdots & \alpha_p E_{p}
\end{array}\right],
$$
where the matrices $E_1,E_2,\ldots,E_p$ are diagonal with
determinant equal to $1$ and each is a product of matrices
$D_{k-1}(i,\theta)$ for $\theta\in\Gamma_{k-1}$. The same
permutational similarity brings $P_k^p$ to
$$
\left[\begin{array}{cccc}
  P_{k-1} &0 & \cdots &0 \\
  0 & P_{k-1} & \cdots & 0 \\
  \vdots & \vdots & \ddots & \vdots \\
  0 & 0 & \cdots & P_{k-1}
\end{array}\right].
$$
Now $P_{k-1}$ and the diagonal blocks generate a metabelian group
in $GL_{p^{k-1}}(F)$ generated by a big cycle and diagonal matrices that are scalar multiples of matrices with the determinant equal to $1$. Recall that property \s\ depends only on $2$-generated subgroups by Proposition \ref{p3.1} and does not depend on multiplication of elements of the group by scalars. Hence, the remaining case then follows by induction.
\hfill\qed\vskip 5pt

\begin{lem}\label{meta4}
Suppose that $\mathcal{G}\subseteq GL_{p^k}(F)$ is an irreducible monomial
$p$-group of class at most $p-1$ that is generated by a
big cycle and a diagonal matrix. Then it has property \s.
\end{lem}

{\em Proof}. By Proposition \ref{bigcycle} it follows that the
group
$$
\widetilde{\mathcal{G}}=\left\{\theta A;\ A\in\mathcal{G}, \theta
\in F,\ \det (\theta A)=1\right\}
$$
has property \s. The lemma now follows since property \s\ does
not depend on multiplication of each element of the group by
scalars.\hfill\qed\vskip 5pt

Assume that $c$ and $e$ are positive integers. Suppose further
that $A$ is the direct product of $c$ copies of the cyclic group
of order $p^e$, and that $\{a_1,a_2,\ldots,a_c\}$ is a set of
generators of $A$. We denote by $B_p(c,e)$ the split extension of
$A$ by an automorphism of order $p^e$ defined by relations:
$b^{-1}a_ib=a_ia_{i+1}$, $i=1,2,\ldots,c-1,$ and
$b^{-1}a_cb=a_c$. It is easy to see that $B_p(c,e)$ is a
metabelian group of exponent $p^e$ and class $c$ and that it is
generated by $a=a_1$ and $b$. The groups $B_p(c,e)$ are called
{\em basic} groups.

Weichsel \cite[p. 62]{Weich1} (see also Brisley \cite{Bris})
showed that each finite me\-ta\-bel\-ian $p$-group of class at
most $p-1$ is in the variety generated by a finite number of the
basic groups $B_p(c,e)$, $c\le p-1$.

\begin{prop}\label{basic}
Suppose that $G=B_p(c,e)$ is a basic metabelian $p$-group with $c\le p-1$ and
that
$\mathcal{G}\subset GL_{p^k}(F)$ is an irreducible representation
of $G$. Then $\mathcal{G}$ is a $2$-generated monomial $p$-group
such that one of the generators is a big cycle and the other is a
diagonal matrix.
\end{prop}

{\em Proof}. Assume that $\psi:G\to GL_{p^k}(F)$ is an irreducible
representation and that $k\ge 1$. Denote by $\mathcal{G}$ the
image of $\psi$. Since $G$ is generated by two elements it
follows that $\mathcal{G}$ is also generated by two elements.
Observe that $A$ is an abelian normal subgroup of $G$ of index
$p^e$. By \cite[Prop. 24, p. 61]{Serre} and arguments in the
proof of \cite[Thm. 16, pp. 66-67]{Serre} it follows that $\psi$
is induced from a representation of $A$. Since $b^iA$,
$i=0,1,\ldots,p^e-1$, is a complete set of cosets of $A$ and
since $\rho$ is irreducible, $\mathcal{G}$ is monomial and one of the generators is diagonal, belonging to $\rho(A)$. The group $G=B_p(c,e)$ is a semi-direct product of $A$ by a cyclic group $C_{p^e}$ of order $p^e$. By \cite[Prop. 25, p. 62]{Serre} all the representations of $G$ are of the type $\theta_{i,\rho}=\chi_i\otimes\rho$, where $\chi_i$ is a representation of $A$, and thus of degree $1$, and $\rho$ a representation of $C_{p^e}$.  Since $k\ge 1$ and $\mathcal{G}$ is irreducible monomial, the image $\psi(b)=\chi_i(1)\otimes\rho(b)$ of the generator $b$ of $C_{p^e}$ is a big cycle.  \hfill\qed\vskip 5pt

Next we prove the main result of the section. First, we introduce
some notation. For two elements $x,y\in G$ we define commutators
$[x,ky]$ inductively as follows: $[x,1y]=[x,y]$ and
$[x,ky]=[[x,(k-1)y],y]$ for $k=2,3,\ldots$

\begin{thm}\label{metabel}
Suppose that $G$ is a metabelian $p$-group. Then the following
are
equivalent:
\begin{enumerate}
  \item $G$ has property \shat,
  \item $G$ is V-regular,
  \item every two generated subgroup of $G$ has class at most
  $p-1$,
  \item the variety of $G$ is generated by a finite group of
  exponent $p$ and a finite group of class at most $p-1$,
  \item $G$ is a $(p-1)$-Engel group, i.e. $[x,(p-1)y]=1$ for and $x,y\in G$,
  \item the variety of $G$ does not contain
  the wreath product of two cyclic groups of order $p$.
\end{enumerate}
\end{thm}

{\em Proof}. The equivalence of (2), (3) and (4) was proved by
Weichsel \cite[Thm. 1.4]{Weich2}. The implication
(1)$\Rightarrow$(2) follows from Theorem \ref{t6.01}.

To prove the implication (3)$\Rightarrow$(1) we may without loss
assume that $G$ is a basic metabelian group $B_p(c,e)$ of class $c\le p-1$. Suppose next that
$\mathcal{G}\subseteq GL_{p^k}(F)$ is an irreducible
representation of $G$. By Proposition \ref{basic}, $\mathcal{G}$ is monomial, generated by $2$ elements
one of which is a big cycle and the other a diagonal matrix. By Lemma
\ref{meta4} it follows that $\mathcal{G}$ has property \s.\

We use \cite[Thm. 3.7]{Groves1} to show that (2) and (3) imply
(5) and (6). Finally \cite[Lem. 3.2, Thms. 3.6 and 3.7]{Groves1}
imply that either (5) or (6) imply (2). \hfill\qed\vskip 5pt

We remark that, in general, the class of a metabelian $p$-group
with property \shat\ can be larger than $p-1$. See, for instance,
the example given by Gupta and Newman in \cite[(3.2)]{GupNew}.
Corollary \ref{pis3} implies that a $3$-group $G$ has property
\shat\ if and only if any of properties (2)--(6) of Theorem
\ref{metabel} holds for $G$. In particular we have:

\begin{cor}
A $3$-group has property \shat\ if and only if it is V-regular.
\end{cor}

\vskip 12pt
\noindent{\bf Acknowledgement:} The authors wish to thank an anonymous referee for careful reading of the original
manuscript and for helpful suggestions.
\bibliographystyle{amsplain}

\end{document}